\documentclass[a4paper,12pt]{article}
\usepackage{times, url}
\usepackage{amsmath,amssymb,amsthm,amsfonts}

\newtheorem{theorem}{Theorem}[section]
\newtheorem{definition}[theorem]{Definition}
\newtheorem{lemma}[theorem]{Lemma}
\begin{document}
\setcounter{page}{1} 

\begin{center}
{\LARGE \bf  On The Log-Concavity of Polygonal Figurate\\[2mm] Number Sequences}
\vspace{8mm}

{\large \bf Fekadu Tolessa Gedefa}
\vspace{3mm}

Department of Mathematics, Salale University\\ 
Fiche, Ethiopia \\
e-mail: \url{fekadu_tolessa@slu.edu.et}
\vspace{2mm}
\end{center}
\vspace{10mm}

\noindent
{\bf Abstract:} This paper presents the log-concavity of the $m$-gonal figurate number sequences. The author gives and proves the recurrence formula for $m$-gonal figurate number sequences and its corresponding quotient sequences which are found to be bounded. Finally, the author also shows that for $m\geq 3$, the sequence $\big \{S_n(m)\big\}_{n\geq 1}$ of $m$-gonal figurate numbers is a log-concave. 
\vspace{10mm}

\section{Introduction} 

Figurate numbers, as well as a majority of classes of special numbers, have long and rich history. They were introduced in Pythagorean school as an attempt to connect Geometry and Arithmetic \cite{deza2012figurate}. A figurate number is a number that can be represented by regular and discrete geometric pattern of equally spaced points \cite{hartman1976figurate}. It may be, say, a  polygonal, polyhedral or polytopic number if the arrangement form a regular polygon, a regular polyhedron or a regular polytope, respectively. In particular, polygonal numbers generalize numbers which can be arranged as a triangle(triangular numbers), or a square (square numbers), to an $m$-gon for any integer $m\geq3$ \cite{weaver1974figurate}.\\
 Some scholars have been studied the log-concavity(or log-convexity) of different numbers sequences such as Fibonacci and Hyperfibonacci numbers, Lucas and Hyperlucas numbers, Bell numbers,Hyperpell numbers, Motzkin numbers, Fine numbers, Franel numbers of order $3$ and 4, Ap\'{e}ry numbers, Large Schr\"{o}der numbers, Central Delannoy numbers, Catalan-Larcombe-French numbers sequences, and so on. See for instance \cite{ahmiaab2014log,aigner1998motzkin,deutsch2001survey,gessel1982some,santana2006some,sun2011delannoy,sun2013congruences,sloane2003line,stanley1989log}.\\
 To the best of the author's knowledge, among all the aforementioned works on the log-concavity and log-convexity of numbers sequences, no one has studied the log-concavity(or log-convexity) of $m$-gonal figurate number sequences. Hence this paper presents the log-concavity behavior of $m$-gonal figurate number sequences.\\
The paper is structured as follows. Definitions and mathematical formulations of figurate numbers are provided in Section 2. Section 3 focuses on the log-concavity of figurate number sequences, and Section 4 is about the conclusion.

 \section{Definitions and Formulas of Figurate Numbers}
 In \cite{deza2012figurate,hartman1976figurate,weaver1974figurate}, some properties of figurate numbers are given. In this paper we continue discussing the properties of $m$-gonal figurate numbers. Now we recall some definitions involved in this paper.
 \begin{definition}
 Let $\big\{s_n\big\}_{n\geq 0}$ be a sequence of positive numbers. If for all $j\geq 1$,  $s_{j}^2\geq s_{j-1}s_{j+1}(s_{j}^2\leq s_{j-1}s_{j+1})$, the sequence $\big\{s_n\big\}_{n\geq 0}$ is called a \it{log-concave}({\it{or a log-convex}}).
 \end{definition}
  \begin{definition}
    Let $\big\{s_n\big\}_{n\geq 0}$ be a sequence of positive numbers. The sequence $\big\{s_n\big\}_{n\geq 0}$ is log-concave(log-convex) if and only if its quotient sequence $\bigg\{\displaystyle\frac{s_{n+1}}{s_n}\bigg\}_{n\geq 0}$ is non-increasing(non-decreasing).
    \end{definition}
    Log-concavity and log-convexity are important properties of combinatorial sequences and they play a crucial role in many fields for instance economics, probability, mathematical biology, quantum physics and white noise theory \cite{asai2001roles,ahmiaab2014log,asai1999bell,dovslic2005log,liu2007log,wang2007log,zhao2014log,zheng2014log,stanley1989log}.\\
    Now we are going to consider the sets of points forming some geometrical figures on the plane. Starting from a point, add to it two points, so that to obtain an equilateral triangle. Six-points equilateral triangle can be obtained from three-points triangle by adding to it three points; adding to it four points gives ten-points triangle, etc.. Then organizing the points in the form of an equilateral triangle and counting the number of points in each such triangle, one can obtain the numbers $1, 3, 6, 10, 15, 21, 28, 36,45, 55,\cdots $, OEIS(Sloane's A000217), which are called {\it{triangular numbers}}, see \cite{bell2005euler,ono1995representation,olson1983triangular,wunderlich1962certain}. The $n^{th}$ triangular number is given by the formula
    \begin{equation}
    S_n=\frac{n(n+1)}{2},n\geq 1.
    \end{equation}
  Following similar procedure, one can construct {\it {square, pentagonal, hexagonal,heptagonal,octagonal,\\nonagonal,decagonal numbers, $\dots$, $m$-gonal}} numbers if the arrangement forms a regular $m$-gon \cite{deza2012figurate}. The $n^{th}$ term  $m$-gonal number denoted by $S_n(m)$ is the sum of the first $n$ elements of the arithmetic progression
    \begin{equation}
    1,1+(m-2),1+2(m-2),1+3(m-2),\dots ,m \geq 3.\label{Sm1}
        \end{equation}
        \begin{lemma}[\cite{deza2012figurate}]
    For $m\geq 3$ and $n\geq 1$,the $n^{th}$ term of $m$-gonal figurate number is given by
    \begin{equation}
    S_n(m)=\displaystyle \frac{n}{2}\big[(m-2)n-m+4\big].\label{Sm2}
    \end{equation}
    \begin{proof}
    To prove (\ref{Sm2}),it suffices to find the sum of the first $n$ elements of (\ref{Sm1}). Hence the first $n$ elements of the arithmetic progression given in (\ref{Sm1}) is: 
    \begin{equation*}
    1,1+(m-2),1+2(m-2),1+3(m-2),\dots,1+(n-1)(m-2),\forall m \geq 3.       \end{equation*}
    Since the sum of the first $n$ elements of an arithmetic progression $s_1,s_2,s_3,\dots,s_n$ is equal to $\displaystyle\frac{n}{2}\big[s_1+s_n\big]$, it follows that
    \begin{align*}
    S_n(m)&=\displaystyle\frac{n}{2}\big[s_1+s_n\big]\\
    &=\frac{n}{2}\big[1+(1+(n-1)(m-2))\big]\\
    &=\frac{n}{2}\big[2+(m-2)n-m+2\big]\\
    &=\frac{n}{2}\big[(m-2)n-m+4\big]~~~~~~~~~~~or\\
    S_n(m)&=\bigg(\frac{m-2}{2}\bigg)\big[n^2-n\big]+n
    \end{align*}
    This completes the proof.
    \end{proof}
    \end{lemma}
    \begin{lemma}[\cite{deza2012figurate}]\label{lemma1}
    For $m\geq 3$ and $n\geq 1$,the following recurrence formula for $m$-gonal numbers hold:
    \begin{equation}
    S_{n+1}(m)=S_n(m)+(1+(m-2)n),S_1(m)=1.
    \end{equation}
    \begin{proof}
    By definition, we have
    \begin{equation*}
    S_n(m)=1+(1+(m-2))+(1+2(m-2))+\cdots + (1+(m-2)(n-2))+(1+(m-2)(n-1))\\
    \end{equation*}
    It follows that
    \begin{align*}
    S_{n+1}(m)&=1+(1+(m-2))+(1+2(m-2))+\cdots + (1+(m-2)(n-1))+(1+(m-2)n)\\
    S_{n+1}(m)&=S_n(m)+(1+(m-2)n).
    \end{align*}
    Thus,for $m\geq 3$ and $n\geq 1$, $$S_{n+1}(m)=S_n(m)+(1+(m-2)n),S_1(m)=1.$$
    \end{proof}
    \end{lemma}
          
\section{Log-Concavity of $m$-gonal Figurate Number Sequences}
In this section, we state and prove the main results of this paper.
\begin{theorem}
    For $m\geq 3$ and $n\geq 3$,the following recurrence formulas for $m$-gonal number sequences hold:
    \begin{equation}
    S_{n}(m)=R(n)S_{n-1}(m)+T(n)S_{n-2}(m)\label{rr1}
    \end{equation}
    with the initial conditions $S_1(m)=1, S_2(m)=m$ and the recurrence of its quotient sequence is given by
    \begin{equation}
    x_{n-1}= R(n)+\displaystyle\frac{T(n)}{x_{n-2}}\label{rr2}
    \end{equation}
    with the initial conditions $x_1=m$, where
    $$R(n)=\displaystyle\frac{m+2(n-2)(m-2)}{1+(n-2)(m-2)}$$ and $$T(n)=-\displaystyle\frac{m-1+(n-2)(m-2)}{1+(n-2)(m-2)}.$$
    \begin{proof}
    By Lemma \ref{lemma1}, we have
    \begin{equation}
    S_{n+1}(m)=S_n(m)+(1+(m-2)n)\label{rr3}
    \end{equation}
    It follows that
    \begin{align}
     S_{n+2}(m)=S_{n+1}(m)+(m-1+(m-2)n)\label{rr4}
    \end{align}
    Rewriting (\ref{rr3}) and (\ref{rr4}) for $n\geq 3$, we have
                \begin{align}
        S_{n-1}(m)&=S_{n-2}(m)+(1+(m-2)(n-2))\label{rr5}\\
        S_{n}(m)&=S_{n-1}(m)+(m-1+(m-2)(n-2))\label{rr6}
        \end{align}
        Multiplying (\ref{rr5}) by  $m-1+(m-2)(n-2)$ and (\ref{rr6}) by  $1+(m-2)(n-2)$, and subtracting as to cancel the non homogeneous part, one can obtain the homogeneous second-order linear recurrence for $S_n(m)$:
        \begin{equation*}
        S_n(m)=\bigg[\displaystyle\frac{m+2(n-2)(m-2)}{1+(n-2)(m-2)}\bigg]S_{n-1}(m)-\bigg[\displaystyle\frac{m-1+(n-2)(m-2)}{1+(n-2)(m-2)}\bigg]S_{n-2}(m),\forall n,m\geq 3.
        \end{equation*}
        By denoting $$\displaystyle\frac{m+2(n-2)(m-2)}{1+(n-2)(m-2)}=R(n)$$ and $$-\displaystyle\frac{m-1+(n-2)(m-2)}{1+(n-2)(m-2)}=T(n),$$ one can obtain
        \begin{equation}
        S_n(m)=R(n)S_{n-1}(m)+T(n)S_{n-2}(m),\forall n,m\geq 3\label{rr7}
        \end{equation}
        with given initial conditions $S_1(m)=1$ and $S_2(m)=m$.
        
        By dividing (\ref{rr7}) through  by $S_{n-1}(m)$, one can also get the recurrence of its quotient sequence $x_{n-1}$ as
        \begin{equation}
        x_{n-1}=R(n)+\displaystyle\frac{T(n)}{x_{n-2}},n\geq 3\label{rr8}
        \end{equation}
        with initial condition $x_1=m.$
        \end{proof}
    \end{theorem}
    \begin{lemma}
    For $m\geq 3$, the $m$-gonal figurate number sequence $\{S_n(m)\}_{n\geq 1}$, let \\$x_n=\displaystyle\frac{S_{n+1}(m)}{S_n(m)}$ for $n\geq 1$. Then we have $1<x_n\leq m$ for $n\geq 1$.
    \begin{proof}
    
      It is clear that 
   \begin{equation*}
   x_1=m,x_2=3-\displaystyle\frac{3}{m},x_3=2-\displaystyle\frac{2}{3(m-1)}>1,\;  \text{for} \; m\geq 3.
   \end{equation*}
   Assume that $x_n>1$ for all $n\geq 3$. It follows from (\ref{rr8}) that
   \begin{equation}
   x_n=\displaystyle\frac{m+2(n-1)(m-2)}{1+(n-1)(m-2)}-\displaystyle\frac{m-1+(n-1)(m-2)}{(1+(n-1)(m-2))x_{n-1}},n\geq 2\label{rr9}
   \end{equation}
   For $n\geq 3$, by (\ref{rr9}), we have
   \begin{align}
   x_{n+1}-1 &=\displaystyle\frac{m-1+n(m-2)}{1+n(m-2)}-\displaystyle\frac{m-1+n(m-2)}{1+n(m-2))x_{n}}\\
   &=\displaystyle\frac{(m-1+n(m-2))x_n-(n(m-2)+m-1))}{(1+n(m-2))x_{n}}\\
   &=\displaystyle\frac{(m-1+n(m-2))(x_n-1)}{(1+n(m-2))x_{n}}\\
      &>0 \;\text{for}  \; m \geq 3.\nonumber
   \end{align}
   Hence $x_n> 1$ for $n\geq 1$ and $m\geq 3.$
       
   Similarly, it is known that 
   \begin{equation}
   x_1=m,x_2=3-\displaystyle\frac{3}{m},x_3=2-\displaystyle\frac{2}{3(m-1)}<m,\;  \text{for} \; m\geq 3.
   \end{equation}
   Assume that $x_n\leq m$ for all $n\geq 3$. It follows from (\ref{rr8}) that
   \begin{equation}
   x_n=\displaystyle\frac{m+2(n-1)(m-2)}{1+(n-1)(m-2)}-\displaystyle\frac{m-1+(n-1)(m-2)}{(1+(n-1)(m-2))x_{n-1}},n\geq 2\label{rr10}
   \end{equation}
   For $n\geq 3$, by (\ref{rr10}), we have
   \begin{align}
   x_{n+1}-m &=-\displaystyle\frac{n(m-2)^2}{1+n(m-2)}-\displaystyle\frac{m-1+n(m-2)}{1+n(m-2))x_{n}}\\
   &=-\displaystyle\frac{n(m-2)^2x_n+n(m-2)+m-1)}{(1+n(m-2))x_{n}}\\
   &< -\displaystyle\frac{n(m-2)^2+n(m-2)+m-1)}{(1+n(m-2))x_{n}}\\
   &=-\displaystyle\frac{n(m-2)(2m-3)}{(1+n(m-2))x_n}\\
   &<0 \; \;   \text{for }  m \geq 3.\nonumber
   \end{align}
   Hence $x_n\leq m$ for $n\geq 1$ and $m\geq 3.$\\
   Thus, in general, from the above two cases it follows that $1<x_n\leq m$ for $n\geq 1$ and $m\geq 3$.
    \end{proof}
    \end{lemma}
    \begin{lemma}[\cite{dovslic2005log}]\label{lemma2}
    Let $\big \{A_n\big\}_{n\geq 0}$ be a sequence of positive real numbers given by the recurrence \[A_n=R(n)A_{n-1}+T(n)A_{n-2},n\geq 2\] with given initial conditions $A_0, A_1$ and $\big \{x_n\big\}_{n\geq 1}$ its quotient sequence, given by \[x_n=R(n)+\displaystyle\frac{T(n)}{x_{n-1}},n\geq 2\] with initial condition $x_1=\displaystyle\frac{A_1}{A_0}$. If there is $n_0\in \mathbb{N}$ such that $x_{n_0}\geq x_{n_0+1}$, $R(n)\geq 0,T(n)\leq 0$, and \[\Delta R(n)x_{n-1}+\Delta T(n)\leq 0\] for all $n\geq n_0$, then the sequence $\big \{A_n\big\}_{n\geq n_0}$ is a log-concave.
       \end{lemma}   
    \begin{theorem}
For all $m\geq 3$, the sequence $\big \{S_n(m)\big\}_{n\geq 1}$ of $m$-gonal figurate numbers is a log-concave.
\begin{proof}
Let $\big \{S_n(m)\big\}_{n\geq 1}$  be a sequence of $m$-gonal figurate numbers given by the recurrence (\ref{rr1}) and let $\big \{x_n\big\}_{n\geq 1}$ be its quotient sequence given by (\ref{rr2}).\\
In order to prove the log-concavity of $\big \{S_n(m)\big\}_{n\geq 1}$ for all $m\geq 3$, by Lemma \ref{lemma2}, we only need to show that $\big \{x_n\big\}_{n\geq 1}$ is non-increasing, $R(n)\geq 0$, $T(n)\leq 0$, and \[\Delta R(n)x_{n-2}+\Delta T(n)\leq 0\] for all $n\geq 3$ .\\
By (\ref{rr7}), since $R(n)\geq 0$ and $T(n)\leq 0$, for $m,n\geq 3$, assume , inductively that $x_1\geq x_2\geq x_3\geq \cdots \geq x_{n-2}\geq x_{n-1}.$\\
Expressing $x_{n}$ from (\ref{rr2}) and taking in to account that $\displaystyle\frac{T(n+1)}{x_{n-1}}\leq \frac{T(n+1)}{x_{n-2}}$, one can obtain
\begin{equation}
x_{n}=R(n+1)+\displaystyle\frac{T(n+1)}{x_{n-1}}\leq R(n+1)+\frac{T(n+1)}{x_{n-2}}\label{rr11}
\end{equation}
Now, we need to show that $x_{n}\leq x_{n-1}.$ To show this, consider
\begin{equation}
R(n+1)+\frac{T(n+1)}{x_{n-2}}\leq R(n)+\frac{T(n)}{x_{n-2}}=x_{n-1}\label{rr12}
\end{equation}
Hence from (\ref{rr11}) and (\ref{rr12}),we can conclude that the quotient sequence $\big \{x_n\big\}_{n\geq 1}$ is non-increasing. 
 It follows from (\ref{rr12}) that
 \begin{equation}
 \big[R(n+1)-R(n)\big]x_{n-2}+ T(n+1)-T(n)\leq 0\label{rr13}
 \end{equation}
 By denoting $R(n+1)-R(n)=\Delta R(n)$ and $T(n+1)-T(n)=\Delta T(n)$, we get the compact expression for (\ref{rr13}) as:
 \begin{equation*}
\Delta R(n)x_{n-2}+\Delta T(n)\leq 0, \forall n\geq 3.
 \end{equation*}
 Thus, by Lemma \ref{lemma2}, the sequence $\big \{S_n(m)\big\}_{n\geq 1}$ of $m$-gonal figurate numbers is a log-concave for $m\geq 3$.\\ This completes the proof of the theorem.
 \end{proof}
 \end{theorem}
\section{Conclusion}
In this paper, we have discussed the log-behavior of $m$-gonal figurate number sequences. We have also proved that for $m\geq 3$, the sequence $\big \{S_n(m)\big\}_{n\geq 1}$ of $m$-gonal figurate numbers is a log-concave.\\[10mm]

\section*{Acknowledgements} 

The author is grateful to the anonymous referees for their valuable comments and suggestions.
 
\makeatletter
\renewcommand{\@biblabel}[1]{[#1]\hfill}
\makeatother

\end{document}